\magnification=\magstep1
\overfullrule=0pt
\def\G{\Gamma}

\def\n{{\bf n}}
\def\psk{{\prod\limits^{s-1}_{k=1}}}
\def\x{{\bf x}}
\def\m{{\bf m}}
\def\lm{{\lambda}}
\def\ita{{i\theta}}
\def\ta{\theta}
\def\a{\alpha}
\def\ty{\infty}
\centerline{\bf Some Systems of Multivariable Orthogonal}
\centerline{\bf Askey-Wilson Polynomials}
\medskip
\centerline{George Gasper\footnote{*}{Dept. of Mathematics, Northwestern University, Evanston, 
IL 60208} and Mizan Rahman\footnote{\dag}{School of Mathematics and Statistics, Carleton University, Ottawa, ON, K1S 5B6,
CANADA. This work was supported, in part, by an NSERC grant \#A6197.}}
\bigskip
\centerline{\bf Abstract}
\medskip
{\narrower\narrower
In 1991 Tratnik derived two systems of multivariable  orthogonal Wilson polynomials and 
considered their limit cases. $q$-Analogues of these systems are derived, yielding systems 
of multivariable orthogonal Askey-Wilson polynomials and
their special and limit cases.}
\bigskip
\noindent
{\bf 1. Introduction}.\quad In [10] Tratnik extended the Wilson [12] polynomials
$$\leqalignno{
&w_n (x;a,b,c,d) = (a+b)_n(a+c)_n (a+d)_n&(1.1)\cr
&\qquad\qquad\qquad\quad\ \times {}_4F_3\left[\matrix{
-n, n+a+b+c+d-1, a+ix, a-ix\cr
a+b, a+c, a+d\cr} ;1\right]&\cr
}
$$
to the multivariable Wilson polynomials (in a different notation)
$$\leqalignno{
&&(1.2)\cr
W_\n(\x) &= W_\n (\x; a, b, c, d, a_2, a_3, \ldots, a_s)&\cr
&= \left[\psk w_{n_k} (x_k; a + \a_{2,k} + N_{k-1}, b + \a_{2,k} + N_{k-1}, a_{k+1} + ix_{k+1}, a_{k+1} - ix_{k+1})\right]&\cr
&\quad \times w_{n_s} (x_s; a + \a_{2,s} + N_{s-1}, b+\a_{2,s} + N_{s-1}, c,d),\cr
}
$$
where, as elsewhere,
$$\x =(x_1,\ldots, x_s),\ \n = (n_1,\ldots, n_s),\ \a_{j,k} = \sum^k_{i=j} a_i,\ \a_k = \a_{1,k},\leqno (1.3)$$
$$N_{j,k} = \sum^k_{i=j} n_i, \ N_k = N_{1,k},\ \a_{k+1,k} = N_{k+1,k} = 0,\ \ 1 \le j \le k \le  s .$$
These polynomials are of total degree $N_s$ in the variables $y_1, \ldots, y_s$ with $y_k = x^2_k$, 
$k = 1, 2, \ldots, s$, and they form a complete set for polynomials in these variables.
\medskip
In Askey and Wilson [1], [2] the notations $W_n (x^2; a, b,c,d)$ and $p_n(-x^2)$ are used for the 
polynomials in (1.1) and their orthogonality relation is given. Tratnik [10, (2.5)] proved that the $W_{\n} (\x)$ polynomials
satisfy the orthogonality relation
$$\int^\ty_{-\ty} \cdots \int^\ty_{-\ty} W_{\n} (\x) W_{\m} (\x) \rho (\x) \ dx_1 \cdots dx_s = \lm_{\n}\, \delta_{\n,
\m}\leqno (1.4)$$
for Re$(a,b,c,d, a_2, \ldots, a_s)> 0$ with
$$\leqalignno{
&\rho (\x) = \G (a+ix_1) \G (a-ix_1) \G(b + ix_1) \G (b-ix_1)&(1.5)\cr
&\quad\quad\quad\times \left[\psk {\G(a_{k+1} + ix_{k+1} + ix_k) \G (a_{k+1} - ix_{k+1} - ix_k) 
\over \G(2ix_k)}\right. &\cr
&\quad\quad\quad\times \left. { \G (a_{k+1} + ix_{k+1} - ix_k)
\G(a_{k+1} - ix_{k+1} + ix_k)\over \G(-2ix_k)} \vphantom{\psk}  \right]&\cr
&\quad\quad\quad\times {\G (c + ix_s) \G (c-ix_s) \G(d + ix_s) \G(d-ix_s)\over
\G(2ix_s) \G(-2ix_s)},&\cr
}
$$
$$\leqalignno{
\lm_{\n} &= (4\pi)^s \left[ \prod^s_{k=1} n_{k}!\, (N_k + N_{k-1} + 2\a_{k+1} -1)_{n_k}\right. &(1.6)\cr
&\quad \times\left. {\G(N_k + N_{k-1} + 2\a_k) \G(n_k + 2a_{k+1})\over\G(2N_k + 2\a_{k+1})}\vphantom{\psk}  \right]&\cr
&\quad \times \G(a + c + \a_{2,s} + N_s) \G(a+d + \a_{2,s} + N_s) \G(b+c + \a_{2,s} + N_s)&\cr
&\quad \times \G(b+d + \a_{2,s} + N_s),&\cr
}
$$
and $2a_1 = a+b$, $2a_{s+1} = c+d$.
\medskip
Tratnik showed that these polynomials contain multivariable Jacobi, Meixner-Pollaczek, Laguerre, 
continuous Charlier, and Hermite polynomials as limit cases, and he used a permutation of the 
parameters  and variables in (1.2) and
(1.4) to show that the polynomials
$$\leqalignno{
&&(1.7)\cr
\tilde W_{\n} (\x) &= \tilde W_{\n} (\x; a, b, c, d, a_2, a_3, \ldots, a_s)&\cr
&= w_{n_1} (x_1; c+ \a_{2,s} + N_{2,s}, d+ \a_{2,s} + N_{2,s}, a,b)&\cr
&\quad\times  \prod^s_{k=2} w_{n_k} (x_k; c + \a_{k+1,s} + N_{k+1,s}, 
d + \a_{k+1,s} + N_{k+1,s}, a_k + ix_{k-1}, a_k - ix_{k-1}) &\cr
}
$$
also form a complete system of multivariable polynomials of total degree $N_s$ in the 
variables $y_k = x^2_k$, $k = 1, \ldots, s$, that is orthogonal with respect to the weight 
function $\rho (\x)$ in (1.5), and with the normalization
constant
$$\leqalignno{
&\tilde\lm_{\n} = (4\pi)^s \Bigg[\prod^s_{k=1} n_{k}! (N_{k,s} + N_{k+1,s} + 2\a_{k,s+1} -1)_{n_k}&(1.8)\cr
&\quad\quad\times {\G(N_{k,s} + N_{k+1,s} + 2\a_{k+1,s+1}) \G(n_k + 2 a_k)\over
\G (2N_{k,s} + 2\a_{k,s+1})}\Bigg]&\cr
&\quad\quad \times \G (a+c + \a_{2,s} + N_s) \G(a+d + \a_{2,s} + N_s) &\cr
&\quad\quad \times \G(b+c+\a_{2,s} + N_s) \G (b+d + \a_{2,s} + N_s).&\cr
}
$$
\indent
The Askey-Wilson polynomials defined as in [1] and [3] by
$$\leqalignno{
p_n (x|q) &= p_n(x;a,b,c,d|q)&(1.9)\cr
&= a^{-n} (ab, ac, ad;q)_n \ _4\phi_3 \left[\matrix{
q^{-n}, abcdq^{n-1}, ae^{\ita}, ae^{-\ita}\cr
ab, ac, ad\cr} ;q,q\right],&\cr
}
$$
where $x = \cos\ta$, are a $q$-analogue of the Wilson polynomials (for the definition of 
the $q$-shifted factorials and the basic hypergeometric series $_4\phi_3$ see [3]). These 
polynomials satisfy the orthogonality relation
$$\int^1_{-1} p_n(x|q) p_m (x|q) \rho (x|q) dx = \lm_n(q)\delta_{n,m} \leqno (1.10)$$
with $\max (|q|, |a|, |b|, |c|, |d|)< 1$,
$$\leqalignno{
\rho (x|q) &= \rho (x;a,b,c,d|q)&(1.11)\cr
&= {(e^{2\ita}, e^{-2\ita};q)_\ty (1-x^2)^{-1/2}\over
(ae^{\ita}, ae^{-\ita}, be^{\ita}, be^{-\ita}, ce^{\ita}, ce^{-\ita}, de^{\ita}, de^{-\ita};q)_\ty}&\cr
}
$$
and
$$\leqalignno{
\lm_n(q) &= \lm_n (a,b,c,d|q)&(1.12)\cr
&= {2\pi (abcd;q)_\ty\over
(q, ab, ac, ad, bc, bd, cd;q)_\ty}&\cr
&\quad \times {(q, ab, ac, ad, bc, bd, cd;q)_n (1-abcdq^{-1})\over
(abcdq^{-1};q)_n (1-abcdq^{2n-1})}.&\cr
}
$$
\indent
In this paper we extend Tratnik's systems of multivariable Wilson polynomials to systems of 
multivariable Askey-Wilson polynomials and consider their special cases. Some $q$-extensions 
of Tratnik's [9] multivariable biorthogonal generalization
of the Wilson polynomials are considered in this Proceedings [4]. $q$-Extensions of Tratnik's [11] 
system of multivariable orthogonal Racah polynomials
and their special cases will be considered in a subsequent paper.
\bigskip
\noindent
{\bf 2. Multivariable Askey-Wilson polynomials}.\quad In terms of the Askey-Wilson polynomials 
a $q$-analogue of the multivariable Wilson polynomials can be defined by
$$\leqalignno{
P_{\n} (\x|q) &= P_{\n} (\x; a,b,c,d, a_2, a_3, \ldots, a_s|q)&(2.1)\cr
&= \left[\prod^{s-1}_{k=1} p_{n_k} (x_k; aA_{2,k}q^{N_{k-1}}, bA_{2,k}q^{N_{k-1}},a_{k+1}e^{i\ta_{k+1}}, a_{k+1}e^{-i\ta_{k+1}}|q)\right]\cr
&\quad \times p_{n_s} (x_s; aA_{2,s}q^{N_{s-1}}, bA_{2,s}q^{N_{s-1}}, c, d\,|q)}
$$
where $x_k = \cos\ta_k$, $A_{j,k} = \prod\limits^k_{i=j}a_i, \, A_{k+1,k} = 1$, 
$A_k = A_{1,k},\ 1 \le j \le k \le  s.$  Our main aim in this section is to show that these polynomials satisfy the orthogonality
relation
$$\int^1_{-1}\cdots \int^1_{-1} P_{\n} (\x|q) P_{\m} (\x|q) \rho (\x|q) dx_1 \cdots dx_s = \lm_{\n} (q)\,
\delta_{\n,\m}\leqno (2.2)$$
with $\max(|q|, |a|, |b|, |c|, |d|, |a_2|, \ldots, |a_s|) < 1$,
$$\leqalignno{
&&(2.3)\cr
\rho (\x|q) &= \rho (\x; a,b,c,d,a_2, a_3,\ldots, a_s|q)&\cr
&= (ae^{i\ta_1}, ae^{-i\ta_1}, be^{i\ta_1}, be^{-i\ta_1};q)_\ty^{-1}&\cr
&\quad \times \left[\psk {(e^{2i\ta_k}, e^{-2i\ta_k};q)_\ty (1-x^2_k)^{-1/2}\over
(a_{k+1}e^{i\ta_{k+1}+ i\ta_k}, a_{k+1}e^{i\ta_{k+1}-i\ta_k}, a_{k+1} e^{i\ta_k-i\ta_{k+1}}, 
a_{k+1}e^{-i\ta_{k+1}-i\ta_k};q)_\ty}\right]&\cr
&\quad \times {(e^{2i\ta_s}, e^{-2i\ta_s};q)_\ty (1-x^2_s)^{-1/2}\over
(ce^{i\ta_s}, ce^{-i\ta_s}, de^{i\ta_s}, de^{-i\ta_s};q)_\ty},&\cr
&&(2.4)\cr
\lm_{\n} (q) &= \lm_{\n} (a,b,c,d, a_2, a_3, \ldots, a_s|q)\cr
&= (2\pi)^s\left[\prod^s_{k=1} {(q, A^2_{k+1} q^{N_k + N_{k-1}-1};q)_{n_k} (A^2_{k+1}q^{2N_k};q)_\ty\over
(q, A^2_k q^{N_k + N_{k-1}}, a^2_{k+1}q^{n_k};q)_\ty}\right]&\cr
&\quad \times (acA_{2,s} q^{N_s}, adA_{2,s}q^{N_s}, bcA_{2,s} q^{N_s}, bd A_{2,s}q^{N_s};q)^{-1}_\ty,&\cr
}
$$
where $a^2_1 = ab$ and $a^2_{s+1} = cd.$  The two-dimensional case was considered by
Koelink and Van der Jeugt [6], but they did not give the value of the norm.
 First observe that by
(1.10)---(1.12) the  integration over $x_1$ in (2.2) can be evaluated to obtain that
$$\leqalignno{
&\int^1_{-1} p_{n_1} (x_1; a, b, a_2e^{i\ta_2}, a_2e^{-i\ta_2}|q)p_{m_1} (x_1;a, b, a_2e^{i\ta_2}, a_2e^{-i\ta_2}|q)&(2.5)\cr
&\quad \times \rho (x_1; a, b, a_2e^{i\ta_2}, a_2e^{-i\ta_2}|q) dx_1&\cr
&= \delta_{n_1,m_1} {2\pi (q, aba^2_2 q^{n_1-1};q)_{n_1} (aba^2_2q^{2n_1};q)_\ty\over
(q, abq^{n_1}, a^2_2 q^{n_1};q)_\ty}&\cr
&\quad \times (aa_2q^{n_1} e^{i\ta_2}, aa_2q^{n_1}e^{-i\ta_2}, ba_2q^{n_1}e^{i\ta_2}, ba_2q^{n_1} e^{-i\ta_2};q)^{-1}_\ty.&\cr
}
$$
After doing the integrations over $x_1, x_2, \ldots, x_j$ for a few $j$ one is led to conjecture that
$$\leqalignno{
&&(2.6)\cr
&\int^1_{-1}\cdots\int^1_{-1} P_{\n}^{(j)} (\x|q) P_{\m}^{(j)} (\x|q) \rho^{(j)} (\x|q) dx_1\cdots dx_j&\cr
&= (2\pi)^j \left[ \prod^j_{k=1} \delta_{n_k, m_k} {(q, A^2_{k+1} q^{N_k + N_{k-1}-1};q)_{n_k} (A^2_{k+1}q^{2N_k};q)_\ty\over
(q, A^2_kq^{N_k + N_{k-1}}, a^2_{k+1} q^{n_k} ;q)_\ty}\right]&\cr
&\quad \times \left( aA_{2,j+1} q^{N_j} e^{i\ta_{j+1}}, aA_{2,j+1}q^{N_j} e^{-i\ta_{j+1}}, bA_{2,j+1}q^{N_j}e^{i\ta_{j+1}},
bA_{2,j+1} q^{N_j}e^{-i\ta_{j+1}};q\right)^{-1}_\ty,&\cr }
$$
where
$$\eqalignno{
P_{\n}^{(j)} (\x|q) &= \prod^j_{k=1} p_{n_k} (x_k; aA_{2,k}q^{N_{k-1}}, bA_{2,k}q^{N_{k-1}}, a_{k+1}e^{i\ta_{k+1}}, a_{k+1}e^{-i\ta_{k+1}}|q),\cr
\rho^{(j)}(\x|q) &= (ae^{i\ta_1}, ae^{-i\ta_1}, be^{i\ta_1}, be^{-i\ta_1};q)^{-1}_\ty\cr
&\quad \times \prod^j_{k=1} {(e^{2i\ta_k}, e^{-2i\ta_k};q)_\ty (1-x^2_k)^{-1/2}\over
(a_{k+1}e^{i\ta_{k+1}+ i\ta_k}, a_{k+1} e^{i\ta_{k+1} - i\ta_k}, a_{k+1}e^{i\ta_k - i\ta_{k+1}}, a_{k+1}e^{-i\ta_{k+1} - i\ta_k};q)_\ty}\cr
}
$$
for $j = 1,2, \ldots, s-1.$ To prove this by induction on $j$, suppose that $j < s-1$, 
multiply (2.6) by the $x_{j+1}$-dependent parts of the weight function and orthogonal 
polynomials, and then integrate with respect to $x_{j+1}$ to get
$$\leqalignno{
&&(2.7)\cr
&(2\pi)^j \left[\prod^j_{k=1} \delta_{n_k,m_k} {(q,A^2_{k+1}q^{N_k+N_{k-1}-1};q)_{n_k} (A^2_{k+1}q^{2N_k};q)_\ty\over
(q, A^2_kq^{N_k+N_{k-1}}, a^2_{k+1}q^{n_k};q)_\ty}\right]&\cr
&\quad \times \int^1_{-1} p_{n_{j+1}} (x_{j+1}; aA_{2,j+1} q^{N_j}, bA_{2,j+1}q^{N_j}, a_{j+2}e^{i\ta_{j+2}}, a_{j+2}e^{-i\ta_{j+2}}|q)\cr
&\quad \times p_{m_{j+1}} (x_{j+1}; aA_{2,j+1}q^{N_j}, bA_{2,j+1}q^{N_j}, a_{j+2}e^{i\ta_{j+2}}, a_{j+2}e^{-i\ta_{j+2}}|q)&\cr
&\quad \times \rho(x_{j+1}; aA_{2,j+1} q^{N_j}, bA_{2,j+1}q^{N_j}, a_{j+2}e^{i\ta_{j+2}}, a_{j+2}e^{-i\ta_{j+2}}|q) dx_{j+1}&\cr
&= (2\pi)^{j+1} \left[ \prod^{j+1}_{k=1} \delta_{n_k,m_k} {(q,A^2_{k+1}q^{N_k+ N_{k-1} -1}; q)_{n_k}
(A^2_{k+1}q^{2N_k};q)_\ty\over (q, A^2_k q^{N_k+N_{k-1}},a^2_{k+1}q^{n_k};q)_\ty}\right]&\cr
&\ \ \ \times \left(aA_{2,j+2} q^{N_{j+1}}e^{i\ta_{j+2}}, aA_{2,j+2}q^{N_{j+1}}e^{-i\ta_{j+2}}, bA_{2,j+2}q^{N_{j+1}}e^{i\ta_{j+2}},
bA_{2,j+2}q^{N_{j+1}}e^{-i\ta_{j+2}};q\right)^{-1}_\ty,&\cr }
$$
which is the $j\to j+1$ case of (2.6), completing the induction proof.
\medskip
Now set $j = s - 1$ in (2.6) and use it and (2.5) to find that
$$\leqalignno{
&\int^1_{-1}\cdots \int^1_{-1} P_{\n} (\x|q) P_{\m}(\x|q) \rho (\x|q) dx_1\cdots dx_s&(2.8)\cr
&= (2\pi)^{s-1} \left[\psk \delta_{n_k,m_k} {(q, A^2_{k+1}q^{N_k+N_{k-1}-1};q)_{n_k} (A^2_{k+1}q^{2N_k};q)_\ty\over
(q, A^2_kq^{N_k+N_{k-1}}, a^2_{k+1}q^{n_k};q)_\ty}\right]&\cr
&\quad \times \int^1_{-1} p_{n_s} (x_s; aA_{2,s}q^{N_{s-1}},bA_{2,s}q^{N_{s-1}},c,d|q)&\cr
&\quad \times p_{m_s} (x_s; aA_{2,s}q^{N_{s-1}}, bA_{2,s}q^{N_{s-1}},c,d|q)&\cr
&\quad \times {(e^{2i\ta_s}, e^{-2i\ta_s};q)_\ty(1-x^2_s)^{-1/2}\over
(ce^{i\ta_s}, ce^{-i\ta_s}, de^{i\ta_s}, de^{-i\ta_s};q)_\ty}&\cr
&\quad \times (aA_{2,s}q^{N_{s-1}}e^{i\ta_s}, aA_{2,s}q^{N_{s-1}}e^{-i\ta_s},bA_{2,s}q^{N_{s-1}}e^{i\ta_s},
bA_{2,s}q^{N_{s-1}}e^{-i\ta_s};q)^{-1}_\ty dx_s&\cr &= \lm_{\n} (q)\, \delta_{\n,\m},&\cr
}
$$
where $\lm_{\n} (q)$ is given by (2.4). This completes the proof of (2.2).
\medskip
Note that the integration region and weight function in (2.2) and (2.3) are invariant under the permutation of variables and parameters
$$a\leftrightarrow c,\ b\leftrightarrow d,\ a_{k+1}\leftrightarrow a_{s-k+1},\quad k = 1,2,\ldots, s-1,\leqno (2.9)
$$
$$\ta_k \leftrightarrow \ta_{s-k+1},\quad k = 1,2,\ldots, s.$$
Hence, when these permutations are applied to (2.2) and (2.3) the transformed polynomials 
also form an orthogonal system with the same weight function. Since the 
polynomials $P_{\n}(\x|q)$ in (2.1) are not invariant under (2.9), we obtain a
second system of multivariable orthogonal Askey-Wilson polynomials, which is 
a $q$-analogue of Tratnik's second system (1.7) of multivariable Wilson
polynomials. After doing the permutation $n_k \leftrightarrow n_{s-k+1}$, $k = 1, 
\ldots, s$, the transformed polynomials and the normalization constant
are given by
$$\leqalignno{
&\qquad\tilde P_{\n} (\x|q) = \tilde P_{\n} (\x; a,b,c,d,a_2,a_3, \ldots, a_s|q)&(2.10)\ \ \cr
&\qquad\qquad \ \ \ \ = p_{n_1} (x_1;cA_{2,s}q^{N_{2,s}},dA_{2,s}q^{N_{2,s}}, a, b|q)&\cr
&\qquad\qquad \ \ \ \ \quad \times \left[\prod^s_{k=2} p_{n_k} (x_k; cA_{k+1,s}q^{N_{k+1,s}}, 
dA_{k+1,s}q^{N_{k+1,s}}, ae^{i\ta_{k-1}},
ae^{-i\ta_{k-1}}|q)\right],&\cr 
&\qquad\tilde\lm_{\n}(q)= \tilde\lm_{\n}(a,b,c,d,a_2,a_3,\ldots, a_s|q)&(2.11)\cr
&\qquad\qquad \ \ = (2\pi)^s \left[\prod^s_{k=1} {(q,A^2_{k,s+1}q^{N_{k,s} + N_{k+1,s}-1};q)_{n_k}
(A^2_{k,s+1}q^{2N_{k,s}};q)_\ty\over (q, A^2_{k+1,s+1}q^{N_{k,s} + N_{k+1,s}}, a^2_kq^{n_k};q)_\ty}\right]&\cr
&\qquad\qquad \ \ \quad \times (acA_{2,s}q^{N_s}, adA_{2,s}q^{N_s}, bcA_{2,s}q^{N_s}, bdA_{2,s}q^{N_s};q)^{-1}_\ty,&\cr
}
$$
with $a^2_1 = ab,$ $a^2_{s+1} = cd,$ and $\max(|q|, |a|, |b|, |c|, |d|, |a_2|, |a_3|, \ldots, |a_s|) < 1$. These polynomials
are  of total degree $N_s$ in the variables $x_1, \ldots, x_s$ and they form a complete set.
\medskip
A five-parameter system of multivariable Askey-Wilson polynomials which is associated with a 
root system of type BC was introduced by Koornwinder [7] and studied with four of the parameters generally 
complex in Stokman [8].
\bigskip
\noindent
{\bf 3. Special Cases of (2.2)}.\quad 
First observe that the continuous dual $q$-Hahn polynomial defined by
$$\leqalignno{
&d_n (x;a,b,c|q) = a^{-n} (ab, ac;q)_n \ _3\phi_2\left[\matrix{
q^{-n}, ae^{\ita}, ae^{-i\ta}\cr
ab, ac\cr};q,q\right]&(3.1)\cr
}
$$
is obtained by taking $d=0$ in (1.9) and $x=\cos\ta.$ Since
 $d_n (x;a,b,c|q)$ is symmetric in its parameters by [3, (3.2.3)], we may 
define the multivariable dual $q$-Hahn polynomials by
$$\leqalignno{
D_{\n} (\x|q) &= D_{\n} (\x; a,b, c, a_2, a_3,\cdots, a_s|q)&(3.2)\cr
&= \left[\psk d_{n_k} (x_k; a_{k+1}e^{i\ta_{k+1}}, a_{k+1}e^{-i\ta_{k+1}},aA_{2,k}q^{N_{k-1}}|q)\right]&\cr
&\quad \times d_{n_s} (x_s; b,c,aA_{2,s}q^{N_{s-1}}|q),&\cr
}
$$
with $x_k = \cos\ta_k$ for $k=1,2, \ldots, s$. It follows from the $b=0$ case of 
(2.2)---(2.4) that the orthogonality relation for these polynomials is
$$\int^1_{-1}\cdots\int^1_{-1} D_{\n} (\x|q) D_{\m}(\x|q)\rho (\x|q)dx_1\cdots dx_s = \lm_{\n}(q)\, \delta_{\n,\m}\leqno
(3.3)$$
with
$$\leqalignno{
&&(3.4)\cr
&\rho(\x|q) = \rho (\x;a, b, c, a_2, a_3,\ldots, a_s|q)&\cr
&\qquad \ \ = (ae^{i\ta_1}, ae^{-i\ta_1};q)^{-1}_\ty&\cr
&\qquad \ \ \quad \times\left[\psk {(e^{2i\ta_k}, e^{-2i\ta_k};q)_\ty (1-x_k^2)^{-1/2}\over
(a_{k+1}e^{i\ta_{k+1} + i\ta_k}, a_{k+1}e^{i\ta_{k+1}-i\ta_k}, 
a_{k+1}e^{i\ta_k -i\ta_{k+1}}, a_{k+1}e^{-i\ta_{k+1}-i\ta_k};q)_\ty}\right]&\cr
&\qquad \ \ \quad \times {(e^{2i\ta_s}, e^{-2i\ta_s};q)_\ty (1-x^2_s)^{-1/2}\over
(be^{i\ta_s}, be^{-i\ta_s}, ce^{i\ta_s}, ce^{-i\ta_s};q)_\ty},&\cr
&&(3.5)\cr
&\lm_{\n} (q)= \lm_{\n}(a, b,c,a_2, a_3, \ldots, a_s|q)&\cr
&\qquad \ \ = (2\pi)^s\left[\prod^s_{k=1} (q^{n_k+1}, a^2_{k+1}q^{n_k};q)^{-1}_\ty\right]
 (abA_{2,s}q^{N_s}, acA_{2,s}q^{N_s};q)^{-1}_\ty,&\cr
}
$$
where $a^2_{s+1} = bc$ and $\max (|q|, |a|, |b|, |c|, |a_2|, |a_3|, \ldots, |a_s|)< 1$.
\medskip
By taking the limit $a\to 0$ in (3.2)---(3.5) we can now deduce that
 the multivariable Al-Salam-Chihara polynomials defined by
$$\leqalignno{
S_{\n} (\x|q) &= S_{\n} (\x; b, c, a_2, a_3,\ldots, a_s|q)&(3.6)\cr
&= \left[\psk p_{n_k} (x_k;a_{k+1}e^{i\ta_{k+1}}, a_{k+1}e^{-i\ta_{k+1}}|q)\right]&\cr
&\quad\times p_{n_s} (x_s; b, c|q).&\cr
}
$$
satisfy the orthogonality relation
$$
\int^1_{-1}\cdots \int^1_{-1} S_{\n} (\x|q) S_{\m} (\x|q) \rho (\x|q) dx_1 \cdots dx_s = \lm_{\n} (q)\,
\delta_{\n,\m}\leqno(3.7)$$
with
$$\leqalignno{
&&(3.8)\cr
&\rho (\x|q) = \left[ \psk {(e^{2i\ta_k}, e^{-2i\ta_k};q)_\ty (1-x^2_k)^{-1/2}\over
(a_{k+1}e^{i\ta_{k+1}+ i\ta_k}, a_{k+1}e^{i\ta_{k+1}- i\ta_k}, a_{k+1} e^{i\ta_k - i\ta_{k+1}}, a_{k+1}
e^{-i\ta_{k+1}-i\ta_k};q)_\ty}\right]&\cr &\quad \quad \quad \quad \times {(e^{2i\ta_s}, e^{-2i\ta_s};q)_\ty (1-x^2_s)^{-1/2}\over
(be^{i\ta_s}, be^{-i\ta_s}, ce^{i\ta_s}, ce^{-i\ta_s};q)_\ty},&\cr
&&(3.9)\cr
&\lm_{\n} (q) = (2\pi)^s  \prod^s_{k=1} (q^{n_k+1}, a^2_{k+1}q^{n_k};q)^{-1}_\ty,&\cr
}
$$
where $a^2_{s+1} = bc,$ $\max (|q|,  |b|, |c|, |a_2|, |a_3|, \ldots, |a_s|)< 1,$ and the Al-Salam-Chihara polynomial
$p_n(x;b,c|q)$ is defined by
$$p_n(x;b,c|q) = b^{-n} (bc;q)_n\ _3\phi_2\left[\matrix{
q^{-n}, be^{i\ta}, be^{-i\ta}\cr
bc, 0\cr} ;q,q\right],\leqno (3.10)$$
see [5, (3.8.1)].
\medskip
Setting
$$a = q^{(2\a+1)/4},\ b = q^{(2\a+3)/4},\ c = -q^{(2\beta +1)/4},\ d = -q^{2\beta + 3)/4}\leqno (3.11)$$
in (2.1) and (2.2) gives a multivariable orthogonal extension of the continuous $q$-Jacobi 
polynomials $P_n^{(\a, \beta)}(x|q)$ defined in [3, (7.5.24)], while setting
$$a = q^{1/2}, \ b = q^{\a + 1/2}, \ c= -q^{\beta + 1/2},\ d = -q^{1/2}\leqno (3.12)$$
in (2.1) and (2.2) gives a multivariable orthogonal extension of the $P_n^{(\a,\beta)}(x;q)$ polynomials
defined in (7.5.25). Also, via [3, (7.5.33)] and [3, (7.5.34) with $q\to q^{1/2}$] 
the $\a = \beta = \lm - 1/2$ substitution gives a multivariable
orthogonal extension of the continuous $q$-ultraspherical polynomials $C_n (x;q^\lm|q)$. 
By letting $\lm \to \ty$ when we use (3.12), i.e. set $a = -d
= q^{1/2}$ and $b = c = 0$,  we get a multivariable orthogonal extension of the continuous $q$-Hermite polynomials defined in [3, Ex. 1.28].
\medskip
A multivariable orthogonal extension of the continuous $q$-Hahn polynomials defined by
$$\leqalignno{
&p_n(\cos(\ta + \phi); a, b|q)&(3.13)\cr
&= (a^2, ab, abe^{2i\phi};q)_n (ae^{i\phi})^{-n}\ _4\phi_3\left[\matrix{
q^{-n}, a^2b^2q^{n-1}, ae^{2i\phi + i\ta},ae^{-i\ta}\cr
a^2, ab, abe^{2i\phi}\cr};q,q\right],&\cr
}
$$
see [3, (7.5.43)], is obtained from (2.1)---(2.4) by replacing $a$, $b$, $c$, $d$, $\ta_k$ 
and $x_k = \cos\ta_k$ by $a_1e^{i\phi}$, $a_1e^{-i\phi}$, $a_{s+1}e^{i\phi}$, $a_{s+1}e^{-i\phi}$, $\ta_k + \phi$ and $\cos (\ta_k + \phi)$,
respectively.
\medskip
It is clear that similar special cases of the second system of multivariable orthogonal 
Askey-Wilson polynomials can be obtained by appropriate specialization of the parameters 
in (2.10) and (2.11). Additional systems of multivariable
orthogonal polynomials will be considered elsewhere.
\bigskip
\noindent
{\bf References}
\medskip
\item{1.} R. Askey and J.A. Wilson, A set of orthogonal polynomials that generalize the Racah 
coefficients or $6$-$j$ symbols, SIAM J. Math. Anal. {\bf 10} (1979), 1008--1016.
\medskip
\item{2.} R. Askey and J.A. Wilson,  Some basic hypergeometric orthogonal polynomials that generalize 
Jacobi polynomials, Memoirs Amer. Math. Soc. {\bf 319}, 1985.
\medskip
\item{3.} G. Gasper and M. Rahman, Basic Hypergeometric Series, Cambridge University Press, 1990.
\medskip
\item{4.} G. Gasper and M. Rahman, $q$-Analogues of some multivariable biorthogonal polynomials,   This Proceedings, 2003.
\medskip
\item{5.} R. Koekoek and R.F. Swarttouw, The Askey-scheme of hypergeometric orthogonal 
polynomials and its $q$-analogue, Report 98-17, Delft University of Technology, http://aw.twi.tudelft.nl/\~{} koekoek/research.html, 1998.
\medskip
\item{6.}  H.T. Koelink and J. Van der Jeugt, 
Convolutions for orthogonal polynomials from Lie and quantum algebra
representations,
SIAM J. Math. Anal. 29 (1998), 794--822
\medskip
\item{7.} T.H. Koornwinder, Askey-Wilson polynomials for root systems of type BC, Contemp. Math. {\bf 138} (1992), 189--204.
\medskip
\item{8.} J.V. Stokman, On BC type basic hypergeometric orthogonal polynomials, Trans. Amer. Math. Soc. {\bf 352} (1999),
1527--1579.
\medskip
\item{9.} M.V. Tratnik, Multivariable Wilson polynomials, J. Math. Phys. {\bf 30} (1989), 2001--2011.
\medskip
\item{10.} M.V. Tratnik, Some multivariable orthogonal polynomials of the Askey tableau---continuous families, 
J. Math. Phys. {\bf 32} (1991),
2065--2073.
\medskip
\item{11.} M.V. Tratnik,  Some multivariable orthogonal polynomials of the Askey tableau---discrete families, 
J. Math. Phys. {\bf 32} (1991),
2337--2342.
\medskip
\item{12.} J.A. Wilson, Some hypergeometric orthogonal polynomials, SIAM J. Math. Anal. {\bf 11} (1980), 690--701.
\bye